\documentclass[12pt]{amsart}
\usepackage{amssymb}  
\usepackage{latexsym} 
\usepackage[all]{xy}
\usepackage{comment} 

\newcommand{\Sha}{\mbox{\wncyr Sh}}

\newcommand{\Q}{{\mathbb Q}}

\newcommand{\PP}{{\mathbb P}}
\newcommand{\oh}{{\mathcal O}}
\newcommand{\RR}{{\mathcal R}}
\newcommand{\MS}{{\mathcal S}}
\newcommand{\BA}{{\mathbb A}}
\newcommand{\tensor}{\otimes}

\newcommand{\Kbar}{\bar{K}}
\newcommand{\Magma}{{\sf MAGMA}\relax}
\newcommand{\GL}{\operatorname{GL}}

\newcommand{\Map}{\operatorname{Map}}
\newcommand{\Sym}{\operatorname{Sym}}
\newcommand{\Spec}{\operatorname{Spec}}
\newcommand{\Sel}{\operatorname{Sel}}

\newcommand{\Mat}{\operatorname{Mat}}
\newcommand{\Proj}{\operatorname{Proj}}
\newcommand{\Trd}{\operatorname{Trd}}
\newcommand{\res}{\operatorname{res}}
\newcommand{\Segre}{\operatorname{Segre}}
\newcommand{\pr}{\operatorname{pr}}
\newcommand{\proj}{\operatorname{proj}}
\newcommand{\eps}{\varepsilon}
\newcommand{\Gm}{\mathbb{G}_{\text{\rm m}}}
\newcommand{\isom}{\cong}
\newcommand{\inj}{\hookrightarrow}
\newcommand{\ra}{\longrightarrow}

\newcommand{\Tr}{\operatorname{Tr}}

\newcommand{\divv}{\operatorname{div}}
\newcommand{\Rbar}{\bar{R}}
\newcommand{\Ob}{\operatorname{Ob}}
\newcommand{\Br}{\operatorname{Br}}

\newcommand{\To}{\longrightarrow}
\newcommand{\mult}{\times} 
\newcommand{\Res}{\operatorname{Res}}

\newfont{\wncyr}{wncyr10 at 12pt}
\newfont{\wncyrten}{wncyr10 at 10pt}

\newtheorem{Proposition}{Proposition}[section]
\newtheorem{Theorem}[Proposition]{Theorem}
\newtheorem{Lemma}[Proposition]{Lemma}

\theoremstyle{definition}
\newtheorem{Definition}[Proposition]{Definition}

\begin{document}

\date{20th November 2006}
\title[Explicit $n$-descent on elliptic curves]%
{Explicit $n$-descent on elliptic curves\\II. Geometry}

\author{J.E.~Cremona}
\address{School of Mathematical Sciences,
         University of Nottingham, 
         University Park, Nottingham NG7 2RD, UK}
\email{John.Cremona@nottingham.ac.uk}

\author{T.A.~Fisher}
\address{University of Cambridge,
         DPMMS, Centre for Mathematical Sciences,
         Wilberforce Road, Cambridge CB3 0WB, UK}
\email{T.A.Fisher@dpmms.cam.ac.uk}

\author{C. O'Neil}
\address{Barnard College, Columbia University,
         Department of Mathematics,
         2990 Broadway MC 4418,
         New York, NY 10027-6902,USA}
\email{oneil@math.columbia.edu}

\author{D. Simon}
\address{Universit\'e de Caen, Campus II - Boulevard Mar\'echal Juin,
         BP 5186--14032, Caen, France}
\email{Denis.Simon@math.unicaen.fr}

\author{M. Stoll}
\address{School of Engineering and Science, 
         International University Bremen 
         (Jacobs University Bremen as of spring 2007),
         P.O. Box 750561, 28725 Bremen, Germany}
\email{M.Stoll@iu-bremen.de}

\begin{abstract} 
This is the second in a series of papers in which we study
the $n$-Selmer group of an elliptic curve. 
 In this paper, we show how to realize elements of the $n$-Selmer group 
explicitly as curves of degree $n$ embedded in~$\PP^{n-1}$. The main
tool we use is a comparison between an easily obtained embedding
into~$\PP^{n^2-1}$ and another map into~$\PP^{n^2-1}$ that factors through
the Segre embedding $\PP^{n-1} \times \PP^{n-1} \to \PP^{n^2-1}$.
The comparison relies on an explicit version of the local-to-global
principle for the $n$-torsion of the Brauer group of the base field.
\end{abstract}

\maketitle


\section{Introduction}

This paper is the second in a series of papers discussing
`Explicit $n$-descent on elliptic curves.' Let $E$ be an elliptic curve
over a number field~$K$, and let $n \ge 3$. 
The elements of the $n$-Selmer group of~$E$ may be viewed as 
isomorphism classes of $n$-coverings $C \to E$.
For us, to do an `explicit $n$-descent' means to represent
each such isomorphism class by giving equations for $C$ as a curve
of degree $n$ in $\PP^{n-1}$. 

It is well known that computing the $n$-Selmer group as
an abstract group gives partial information about both the
Mordell-Weil group $E(K)$ and the Tate-Shafarevich group $\Sha(K,E)$. 
There are likewise several motivations for computing 
the $n$-Selmer group in the explicit sense described above.
Firstly, it is known that a point in $E(K)$ of $x$-coordinate 
height~$h$ lifts to a point of height approximately $h/(2n)$ on one 
of the Selmer $n$-coverings
(see Theorem B.3.2 in \cite{Hindry-Silverman}).
Thus explicit $n$-descent enables us to find generators
for the Mordell-Weil group more easily, 
and hence sometimes
to show that the $n$-torsion of the Tate-Shafarevich group is
trivial. Secondly, if we have already computed
the Mordell-Weil group, for instance using descent at some other 
$n'$, then we can use explicit $n$-descent to exhibit 
concrete examples of non-trivial $n$-torsion elements of~$\Sha(K,E)$.
Our work is also likely to form a useful starting point for
performing higher descents, and for computing the Cassels-Tate pairing.

The first paper in this series \cite{paperI} gives more 
background on the theory of descent. In particular 
we explain a number of different interpretations of the elements of the
group $H^1(K, E[n])$, which contains the Selmer group as a subgroup.
One of these interpretations, and the one which is relevant to this paper,
is as `Brauer-Severi diagrams' $[C \to S]$;
such a diagram consists of
a morphism that is a twist of the embedding $E \to \PP^{n-1}$
associated to the complete linear system $|n (\oh)|$, where $\oh$ is
the origin of~$E$. In particular $C$ is a torsor under~$E$ and 
$S$ is a Brauer-Severi variety of dimension~$n-1$. If $C$ has points
everywhere locally, then so does~$S$, hence (by Global Class Field Theory)
$S$ is isomorphic to~$\PP^{n-1}$. Our goal in this paper is to explain
how one can obtain equations for the image of~$C$ in $S \cong \PP^{n-1}$
in this case. 

In brief, we compare two different embeddings
of $C$ into $\PP^{n^2-1}$. The first comes from the fact that, even 
without assuming $C$ has points everywhere locally, the cohomology map 
$H^1(K, E[n]) \to H^1(K, E[n^2])$ sends the Brauer-Severi diagram
$[C \to S]$ to a diagram $[C \to S']$ with $S' \isom \PP^{n^2-1}$.
The second map is more abstract: starting with the diagram $C \to S$ 
we form a dual map $C \to S^\vee$ and thus
$C \to S \times S^\vee.$  We then compose with the (generalised) Segre 
embedding to obtain $$C \to S \times S^\vee \to \PP^{n^2-1}.$$
Equations for the image of the first map to $\PP^{n^2-1}$ 
are given in Section~\ref{CtoPR}.
By a suitable change of coordinates on $\PP^{n^2-1}$
followed by projection to a hyperplane, we obtain equations for
the image of the second map. Finally, we pull back to $C \to S$ and, when
we start with a Selmer group element, obtain equations for 
$C \to \PP^{n-1}$. For obvious reasons we refer to this 
as the Segre embedding method.

We will not be concerned in this paper with the details of implementation;
these will be discussed in the third paper of the series~\cite{paperIII}.
However, we would like to mention that the Segre embedding method, as well
as two more methods discussed in \cite{paperI}, have been implemented
for $n = 3$ and $K=\Q$ and are available as part of the \Magma\ 
computer algebra system \cite{Magma} (version
2.13 or later).  

All three methods rely for their practical implementation on a `Black
Box' that computes, for a given central simple $K$-algebra~$A$ of
dimension~$n^2$ that is known to be isomorphic to the matrix
algebra~$\Mat_n(K)$, an explicit isomorphism $A \to \Mat_n(K)$.
Algorithms for this when $K= \Q$ will be described in~\cite{paperIII}.


\section{Background and overview}\label{bckgd}

Unless stated otherwise, $K$ will denote a number field, with absolute
Galois group~$G_K$, and $E$ will be an elliptic curve over~$K$ with 
origin~$\oh$. Let $n$ be a positive integer. Recall the definition 
of the {\em $n$-Selmer group} of~$E$: the short exact sequence of 
$K$-group schemes
\[ 0 \ra E[n] \ra E \stackrel{[n]}{\ra} E \ra 0 \]
gives rise to the following commutative diagram with exact rows.
\[ \xymatrix{ E(K) \ar[r]^-{\delta} \ar[d]
                & H^1(K, E[n]) \ar[r] \ar[d] \ar[rd]^{\alpha}
                & H^1(K, E) \ar[d] \\
              \prod_v E(K_v) \ar[r]^-{\delta}
                & \prod_v H^1(K_v, E[n]) \ar[r]
                & \prod_v H^1(K_v, E)
            }
\]
Here, $v$ runs through all places of~$K$.
The $n$-Selmer group $\Sel^{(n)}(K, E)$ is then defined to be the
kernel of~$\alpha$.

The following notation and facts can be found in~\cite{paperI}.

We assume that $n \geq 3$. Then there is an 
embedding $f : E \to \PP^{n-1}$ associated to the complete linear system 
$|n (\oh)|$.  
In fact, if $n = 2$ one would work with a double cover, in which case
our algorithm still works with minor changes. Indeed it reduces to 
the classical number field method for 2-descent 
(see for example \cite{CaL}, Lecture 15).

We can view an element of $H^1(K, E[n])$ as a twist of the diagram
$f : E \to \PP^{n-1}$, i.e., as a diagram of the form $[C \to S]$, where 
$C$ is a torsor under~$E$ and $S$ is a Brauer-Severi variety of 
dimension~$n-1$. We call such a diagram a {\em Brauer-Severi diagram.}
In this interpretation, a diagram $[C \to S]$ corresponds to an element of 
the $n$-Selmer group if and only if $C$ has points everywhere locally.

The {\em period-index obstruction map}, defined in \cite{cathy2}, 
is a quadratic map 
\[ \Ob_n : H^1(K, E[n]) \To \Br(K)[n]. \]  It 
sends the diagram $[C \to S]$ to the class of~$S$ in~$\Br(K)[n]$. 
We say an element $\xi$ of~$H^1(K, E[n])$
has {\em trivial obstruction} if $\Ob_n(\xi) = 0$, 
equivalently the corresponding diagram $[C \to S]$ has
$S \cong \PP^{n-1}$.

Alternatively we can view an element $\xi$ of~$H^1(K, E[n])$
as a pair $(C, [D])$ where $[D]$ is a $G_K$-invariant divisor class on~$C.$
The class~$[D]$ is represented by a rational divisor if and only if the
element $\xi$ has trivial obstruction. 

Our algorithm applies not only to elements of the Selmer group,
but more generally to any element in $H^1(K, E[n]) $ with trivial obstruction.

The natural map $H^1(K, E[n]) \rightarrow H^1(K, E[n^2])$ brings
the pair $(C, [D])$ to the pair $(C, [nD]).$  Moreover, composing the above 
with $\Ob_{n^2}$ gives the zero map: there is a commutative diagram
\[ \xymatrix{ H^1(K, E[n]) \ar[rr]^{\Ob_n} \ar[d] & & 
                                   \Br(K)[n] \ar[d]^{\cdot n} \\
              H^1(K, E[n^2]) \ar[rr]^{\Ob_{n^2}} & & \Br(K)[n^2]
            }
\]
In other words, the divisor class $[nD]$
is always represented by a rational divisor of degree~$n^2,$ or equivalently,
the above map takes any diagram $[C \to S]$ to 
a diagram $[C \to \PP^{n^2-1}]$ with trivial obstruction.  It is
relatively easy to find equations for the image of this map. 
We do this in Section~\ref{CtoPR}.

Our basic question then is, starting with an element of $H^1(K, E[n])$
with trivial obstruction, how do we {\it reverse} the map
\[ [C \to \PP^{n-1}] \longmapsto [C \to \PP^{n^2-1}] \;? \]

The diagram $[C \to S]$ naturally extends to give a map
\[ \lambda_C : C \To S \times S^\vee \To \PP(A) \]
where $A$ is the {\em obstruction algebra,} the central simple
$K$-algebra associated to the Brauer group element $\Ob_n(C \to S) =
[S]$ (see~\cite[p.~160]{Serre}). In this paper we describe an
algorithm for writing down both structure constants for~$A$ and
equations for~$C$ as a curve of degree~$n^2$ in $\PP(A) \cong
\PP^{n^2-1}$. In fact we specify the equations in Section~\ref{CtoPR}
and the structure constants in Section~\ref{RcongA}.  

In the case of trivial obstruction, we know that there exists an
isomorphism of $K$-algebras $A \isom \Mat_n(K)$, which is called a
{\em trivialisation} of~$A$. 
Using this isomorphism, we may obtain equations for $C$ in~$\PP^{n-1}$
as a curve of degree~$n$ by projecting onto a column.

\medskip

The algorithm is split into parts, each of which (except the first) correspond 
to a piece of the `master diagram' found in Theorem~\ref{masterdiagram}
below, which we reproduce here.
\[ \xymatrix{ C \ar[r]^-{f'_C \times {f'_C}^\vee} \ar[d]^{g_C}
                & \PP^{n-1} \times (\PP^{n-1})^\vee \ar[r]^-{\Segre} 
                & \PP(\Mat_n) \\
              \PP(R) \ar[r]_-{\varphi_\rho}^-{\sim} 
                & \PP(A_\rho) \ar[r]_-{\tau_\rho}^-{\sim}
                & \PP(\Mat_n) \ar@{-->}[u]^{\proj} }
\]

The first step is to realize $H^1(K, E[n])$ as a subgroup of a concrete
group. Let $R$ be the
affine algebra of the group scheme~$E[n]$. Addition in~$E[n]$ translates
into the comultiplication homomorphism $\Delta : R \to R \tensor_K R$.
We use this to define a group homomorphism
\[ \partial : R^\mult \To (R \tensor_K R)^\mult, \quad
              \alpha \longmapsto \frac{\alpha \tensor \alpha}{\Delta(\alpha)} .
\]
Then there is a natural embedding of~$H^1(K, E[n])$ into
$(R \tensor_K R)^\mult/\partial R^\mult$. See Section~\ref{CtoPR} below
and~\cite[Section~3]{paperI}. We can then compute the Selmer group
as a subgroup of $(R \tensor_K R)^\mult/\partial R^\mult$, see~\cite{paperIII}
for details. Therefore, in the following we can assume that our element
of~$H^1(K, E[n])$ is represented by some $\rho \in (R \tensor_K R)^\mult$.

In the second step, starting with $\rho \in (R \tensor_K R)^\mult$
coming from an element of $H^1(K, E[n])$, we construct an embedding
$g_C : C \to \PP(R)$, where $\PP(R)$ is the projective space associated
to the $K$-vector space~$R$. We need to choose a $K$-basis of~$R$ in order
to write down explicit equations for the image of~$g_C$. 
This step is explained in Section~\ref{CtoPR}.

In the third step, we define a new multiplication on~$R$,
which depends on $\rho$, that turns it
into a central simple $K$-algebra~$A_\rho$. This is the 
obstruction algebra for the element of~$H^1(K, E[n])$ represented by~$\rho$.
In other words, we create an explicit isomorphism of $K$-vector spaces
$\varphi_{\rho}: R \to A_{\rho}$. 
This will be explained in Section~\ref{RcongA}.

The fourth step makes use of the fact that $S \cong \PP^{n-1}$, or
equivalently that $A_{\rho}\cong\Mat_n(K)$, when~$\rho$ represents an
element of~$H^1(K, E[n])$ which has trivial obstruction.  In
Section~\ref{BlackBox}, we give an explicit trivialisation of the
algebra~$A_1$, the central simple algebra coming from the trivial
element of $H^1(K, E[n])$, which serves as a normalisation.  For the
general case a trivialisation map $\tau_\rho : A_\rho \to \Mat_n(K)$
will come from our `Black Box', to be described in more detail in the
third paper in this series~\cite{paperIII}.

In the fifth step, we project the image of $\tau_\rho \circ
\varphi_\rho \circ g_C$ into the hyperplane of trace~zero matrices and
show that our total map factors through the Segre embedding:
\[  
   C \To \PP^{n-1} \times (\PP^{n-1})^\vee \stackrel{\Segre}{\To}
   \PP(\Mat_n) .  
\]
This step makes up Section~\ref{SegreforE} for the case $\rho =1$,
i.e.  $[C\to \PP^{n-1}] \cong [E \to \PP^{n-1}]$.  The general case is
described in Sections~\ref{SegreforC} and~\ref{CinPn}.

The final step of the algorithm is to make use of the Segre
factorisation.  It is here that an explicit trivialisation of the
obstruction algebra is required.  We pull back under the Segre
embedding and project to the first factor, which gives us equations
for $C \to \PP^{n-1}$. This is explained in Section~\ref{CinPn}.


\section{Finding Equations for $C$ in $\PP^{n^2-1}$}\label{CtoPR}

We continue to take $K$ a number field, but in fact
our algorithm works over any field whose characteristic is prime to $n$
(although below, we assume for simplicity that the characteristic is
neither 2 nor~3).
We denote by $[n]$ the multiplication-by-$n$ map on $E$.
We fix a Weierstrass equation 
\[ E : \quad y^2 = x^3 + a x + b \]
and define for $T_1, T_2 \in E[n](\Kbar)$  a rational function
$r_{(T_1,T_2)}$ in $\Kbar(E)^\mult$,
\[ \label{rPweier}
    r_{(T_1,T_2)}(P)
      = \left\{ \begin{array}{ll}
                  1           & \text{if $T_1 = \oh$ or $T_2 = \oh$} \\[1mm]
                  x(P)-x(T_1) & \text{if $T_1 + T_2 = \oh$ and $T_1 \neq \oh$} 
                                \\[1mm]
                  \frac{y(P)+y(T_1+T_2)}{x(P)-x(T_1+T_2)} - \lambda(T_1,T_2)
                              & \text{otherwise,} 
                \end{array}\right.
\]
where $\lambda(T_1,T_2)$ denotes the slope of the line joining
$T_1$ and $T_2$, respectively of the tangent line at~$T_1 = T_2$
if the points are equal.

For $T \in E[n](\Kbar)$, there exists a rational function $G_T \in \Kbar(E)$ with divisor
\[ \divv (G_T) = \sum_{n  S = T} (S) - 
     \sum_{n S=\oh} (S) = [n]^*(T)- [n]^*(\oh) . \]
(See~\cite[Section~III.8]{Silverman}.) 

\begin{Proposition}\label{scaling}
  We can scale the $\{G_T\}$ such that
  \begin{enumerate}\addtolength{\itemsep}{1mm}
    \item The map $T \mapsto G_T$ is $G_K$-equivariant,
    \item For each $T_1,T_2 \in E[n](\Kbar)$ and 
          $P \in E(\Kbar)\backslash E[n^2](\Kbar)$ we have
          \[ r_{(T_1, T_2)}(nP) = \frac{G_{T_1}(P)G_{T_2}(P)}{G_{T_1+T_2}(P)}. 
          \]
    \item $G_\oh =1$, and for $T \neq \oh$ the residue of $G_T$ at $\oh$ with 
          respect to the local parameter~$x/y$ is~$\frac{1}{n}$.
  \end{enumerate}
\end{Proposition}

\begin{proof}  
  Define the scalings of the $G_T$ so that condition~(iii) holds.  That means 
  that  with respect to the local parameter $t= x/y$, when $T \not = \oh$,
  we can write $G_T(t)= \frac{1}{n} t^{-1} + \dots$, where 
  `\dots' signifies `higher order terms.'
  This choice of scaling makes the map $T \mapsto G_T$ visibly 
  $G_K$-equivariant.

  When $T_1= \oh$ or $T_2= \oh,$ condition~(ii) holds trivially.  In any case,
  a calculation of divisors gives
  \[ \divv (r_{(T_1,T_2)}) = (T_1) + (T_2) - (\oh) - (T_1+T_2) \] 
  for all choices of $T_1$ and $T_2$.
  So the divisor of $r_{(T_1, T_2)}\circ [n]$ is
  \[ \divv (r_{(T_1,T_2)}\circ [n])
        = \sum_{nS=T_1} (S) + \sum_{nS=T_2} (S)
            - \sum_{nS=\oh} (S) - \sum_{nS=T_1+T_2} (S) .
  \]  
  This is exactly the divisor of 
  $P \mapsto G_{T_1}(P)G_{T_2}(P)/G_{T_1+T_2}(P)$.
  Therefore condition~(ii) holds up to a scalar. In the following, we consider
  the case $T_1 + T_2 \neq \oh$.
  We write $x(t) = t^{-2} + \dots$ and $y(t) = t^{-3} + \dots$.
  Then locally at~$\oh$,
  \begin{align*}
   r_{(T_1,T_2)}(t)
     &= \frac{y(t)+y(T_1+T_2)}{x(t)-x(T_1+T_2)} - \lambda(T_1,T_2) \\
     &= \frac{t^{-3} + \dots}{t^{-2} + \dots} - \lambda(T_1,T_2) \\
     &= t^{-1} + \dots
  \end{align*}
  Next, from~\cite[Prop.~IV.2.3]{Silverman}, we have $[n](t)= nt + \dots$, 
  hence
  \[ r_{(T_1,T_2)}\circ[n] (t)
       = (nt)^{-1}+ \dots = \frac{1}{n} t^{-1} + \dots
  \]
  Comparing this with
  \[ \frac{G_{T_1}(t)G_{T_2}(t)}{G_{T_1+T_2}(t)}
        = \frac{(\frac{1}{n} t^{-1} + \dots)(\frac{1}{n} t^{-1} + \dots)}%
               {\frac{1}{n} t^{-1} + \dots}
        = \frac{1}{n} t^{-1} + \dots
  \]
  shows that the scalar is~$1$.
  The case $T_1+T_2=\oh$ is similar.
\end{proof}

In preparation for defining the embedding of $C$ in $\PP^{n^2-1}$, 
we recall some facts from~\cite{paperI}.
Let $R$ be the affine algebra of $E[n]$, i.e.,
\begin{equation*} 
  R = \Map_K(E[n](\Kbar), \Kbar). 
\end{equation*}
It is isomorphic to a product of (finite) field extensions of~$K$, one
for each $G_K$-orbit in~$E[n](\Kbar)$. We also work with the algebra
\[ \Rbar = R \tensor_K \Kbar = \Map(E[n](\Kbar),\Kbar). \]

The Weil pairing $e_n: E[n] \times E[n] \to \mu_n$ determines an injection 
\[ w : E[n](\Kbar) \inj \Rbar^\mult = \Map(E[n](\Kbar),\Kbar^\mult) \]
via $w(S)(T) = e_n(S,T)$. 

As in Section \ref{bckgd}, we define $\partial: \Rbar^\mult \to (\Rbar \tensor_{\Kbar} \Rbar)^\mult$ via
\begin{equation} \label{defpartial}
  \partial \alpha = \frac{\alpha \tensor \alpha}{\Delta(\alpha)},
    \qquad\text{i.e.,}\qquad
    (\partial \alpha)(T_1,T_2)
    = \frac{\alpha(T_1) \alpha(T_2)}{\alpha(T_1+T_2)};
\end{equation}
then there is an exact sequence (cf.~\cite[Section~3]{paperI})
\begin{equation} \label{exseq1}
  0 \ra E[n](\Kbar) \stackrel{w}{\ra} \Rbar^\mult 
          \stackrel{\partial}{\ra} (\Rbar \tensor_{\Kbar} \Rbar)^\mult.
\end{equation}

For $V$ a vector space over~$K$, 
we write $\PP(V) = \Proj(K[V])$, where $K[V]= \oplus_{d \ge 0} \Sym^d (V^*)$
is the ring of polynomial functions on $V$.

We define $\RR = \Res_{R/K}(\BA^1)$, or equivalently, $\RR= \Spec(K[R])$.  
For any $K$-scheme $X,$ we have
\[ \RR(X) = \BA^1(\Spec(R) \times_{\Spec(K)} X). \]  
In particular, $\RR(L) = R\otimes_K L$ for any field extension $L/K.$

We also define $\RR^\mult = \Res_{R/K}(\Gm)$ and 
$\MS^\mult = \Res_{R \tensor_K R/K}(\Gm)$.
These schemes inherit a multiplication from~$\Gm$.
The groups of $K$-rational points are $\RR^\mult(K) = R^\mult$
and $\MS^\mult(K) = (R \tensor_K R)^\mult$.
We may identify $\RR^\mult$ with an open subscheme of~$\RR$. 

With this notation, the exact sequence of $G_K$-modules~(\ref{exseq1}),
becomes an exact sequence of $K$-group schemes:
\begin{equation} \label{exseq2}
  0 \ra E[n] \stackrel{w}{\ra} \RR^\mult \stackrel{\partial}{\ra} \MS^\mult .
\end{equation}

Part (i) of Proposition~\ref{scaling} allows us to package the functions~$G_T$ 
to form a scheme map
$g_E: E \to \PP(R)$ sending $P \in E(\Kbar)$ to the class of the map
$T \mapsto G_T(P)$.  Away from the subscheme $E[n^2]$, the map $g_E$ can be lifted to a map $g^0_E$ to $\RR^\mult$. Then we have a commutative diagram:
\[ \xymatrix{ E \setminus E[n^2] \ar@{^(->}[d] \ar[r]^-{g^0_E}
                & \RR^\mult \ar[d] \\
              E \ar[r]^-{g_E} & \PP(R) } 
\]

Next, we use the $G_K$-equivariance of $(T_1, T_2) \mapsto r_{(T_1, T_2)}$
to package the functions $r_{(T_1, T_2)}$ to form a scheme map
$r: E \setminus E[n] \to \MS^\mult$ sending 
$P \in E(\Kbar) \setminus E[n](\Kbar)$ to the map
$r(P): (T_1, T_2) \mapsto  r_{(T_1, T_2)}(P)$.

\begin{Proposition}\label{compose}
  The following diagram commutes.
  \[ \xymatrix{ E \setminus E[n^2] \ar[r]^{[n]} \ar[d]_{g_E^0} 
                 & E \setminus E[n] \ar[d]^{r} \\
                \RR^\mult \ar[r]^{\partial} & \MS^\mult
              }
  \]
\end{Proposition}

\begin{proof}
  We take a geometric point $P \in E(\Kbar) \setminus E[n^2](\Kbar)$
  and observe that 
  \[ \partial(g_E^0(P))(T_1, T_2)
         = \frac{G_{T_1}(P) G_{T_2}(P)}{G_{T_1+T_2}(P)}.
  \]
  By Proposition~\ref{scaling}(ii), this equals $r_{(T_1, T_2)}(nP)$.
\end{proof}

For $T \in E[n](\Kbar)$ we denote by $z_T$ the coordinate function on 
$\RR \times_{\Spec(K)} \Spec(K(T))$ given by evaluating at $T$, so
$z_T(\alpha) = \alpha(T)$.

\begin{Proposition} \label{eqs}
  Given a Weierstrass equation for~$E$, we can explicitly compute a set
  of $n^2 (n^2-3)/2$ linearly independent quadrics over~$K$ which
  define the image of
  \[ g_E : E \ra \PP(R) \cong \PP^{n^2-1} . \]
  If $E[n](\Kbar) = E[n](K)$, then the $z_T$ are coordinate functions on~$\RR$,
  and the defining quadrics can be split into two groups as follows.
  For all $T_1, T_2 \in E[n](\Kbar) \setminus \{\oh\}$, we have 
  \[ \bigl(x(T_1)-x(T_2)\bigr)  z_{\oh}^2
       +  z_{T_1}z_{-T_1} - z_{T_2}z_{-T_2}, 
  \] 
  and for all $T_{11}, T_{12}, T_{21}, T_{22} \in E[n](\Kbar) \setminus \{\oh\}$ 
  such that 
  \[ T_{11} + T_{12} = T_{21} + T_{22} = T \neq \oh, \] 
  we have
  \[ \bigl(\lambda(T_{21}, T_{22})-\lambda(T_{11}, T_{12})\bigr) z_{\oh} z_T  
       - z_{T_{11}}z_{T_{12}} + z_{T_{21}}z_{T_{22}}.
  \]
\end{Proposition}

\begin{proof}
  We first note that the $G_T$ are linearly independent. This follows from
the fact they are eigenfunctions for distinct characters with respect
to the action of $E[n]$ by translation. (We are using  
the definition of
  the Weil pairing in \cite[Section~III.8]{Silverman}.) 
  Since there are $n^2$ functions $G_T$,
  they form a basis for the
  Riemann-Roch space of the divisor $[n]^*(\oh)$.  
  Hence $g_E$ embeds $E$ into $\PP(R) \cong \PP^{n^2-1}$
  as an elliptic normal curve of degree~$n^2$.
  
  Now let $P \in E(\Kbar) \setminus E[n^2](\Kbar)$, 
  and let $z = g_E^0(P) \in \RR^\mult(\Kbar)$ be 
  projective coordinates for~$g_E(P)$.
  By Proposition~\ref{compose}, we then have $r(nP) = \partial z$, 
  or equivalently, $r(nP) \Delta(z) = z \tensor z$.
  We wish to eliminate $P$ from this equation. Since $z_{\oh}(g_E^0(P)) = 1$,
  we can make the equation homogeneous by multiplying the left hand side
  with~$z(\oh)$. 
  This gives
  \[ r(nP) \, z(\oh) \Delta(z) = z \tensor z . \]
  Writing everything out in terms of a $K$-basis of~$R$, we obtain 
  $n^4$~quadrics,
  some of whose coefficients involve rational functions of~$nP$.
  We can eliminate these rational functions by linear algebra over~$K$, 
  to obtain a set of quadrics in $K[R]$ 
  --- this is what we do in practice. 
  In order to determine the dimension of the space spanned by them
  and the geometry
  of the object defined by them, we can work over~$\Kbar$. In fact,
  it is sufficient to work over $L = K(E[n])$.
  
  Over~$L$, the coordinate functions~$z_T$ are defined. In terms of these,
  our system of equations is
  \[ r_{(T_1, T_2)}(nP) \, z_{\oh} z_{T_1 + T_2} = z_{T_1} z_{T_2} , \]
  parametrised by $(T_1, T_2) \in E[n](\Kbar) \times E[n](\Kbar)$.
  
  If $T_1 = \oh$ or $T_2 = \oh$, this reduces to a tautology.
  
  If $T_1 = T \neq \oh$ and $T_2 = -T$, then we get
  \[ \bigl(x(nP) - x(T)\bigr) z_{\oh}^2 = z_{T} z_{-T} . \]
  We can eliminate $x(nP)$ by taking differences. Taking into account
  the symmetry $T \leftrightarrow -T$, this gives us $d_1$ independent
  quadrics, where
  \[ d_1 = \#\Bigl(\frac{E[n](\Kbar) \setminus \{\oh\}}{\{\pm 1\}}\Bigr) - 1
         = \begin{cases}
             (n^2-3)/2 & \text{if $n$ is odd} \\
             n^2/2     & \text{if $n$ is even.}
           \end{cases}
  \]
  
  If $T_1 + T_2 = T \neq \oh$ and $T_1, T_2 \neq \oh$, we obtain
  \[ \Bigl(\frac{y(nP) + y(T)}{x(nP) - x(T)} - \lambda(T_1, T_2)\Bigr)
       z_{\oh} z_T
         = z_{T_1} z_{T_2} .
  \]
  Fixing $T$, we can again eliminate the dependence on~$P$ by
  taking differences. Taking into account the symmetry 
  $(T_1, T_2) \leftrightarrow (T_2, T_1)$, this provides us with $d_2$
  independent quadrics, where
  \begin{align*}
    d_2 &= \#\Bigl(\frac{\{(T_1, T_2) : T_1, T_2, T_1+T_2 \neq \oh\}}%
                        {(T_1, T_2) \sim (T_2, T_1)}\Bigr)
            - \#\bigl(E[n](\Kbar) \setminus \{\oh\}\bigr) \\
        &= \begin{cases}
             (n^2-1)(n^2-3)/2 & \text{if $n$ is odd} \\
             n^2(n^2-4)/2     & \text{if $n$ is even.}
           \end{cases}
  \end{align*}
  Together, we
  obtain $d_1 + d_2 = n^2(n^2-3)/2$ independent quadrics in either case.
  
  We have found an $n^2(n^2-3)/2$-dimensional space of quadrics
  vanishing on the image of~$g_E$.  In general 
  (see for example the Corollary to Theorem 8 in \cite{Mumford}) 
  the homogeneous ideal of an elliptic normal curve of
  degree $m \geq 4$ is generated by a vector space of quadrics of dimension
  $m (m-3)/2$. Therefore, our quadrics define the image of~$E$ in~$\PP(R)$
  under~$g_E$.
\end{proof}

By Section~3 of~\cite{paperI}, we can identify $H^1(K, E[n])$ with a subgroup 
of $(R \tensor_K R)^\mult/ \partial R^\mult$. 
Specifically, we define an injective map
\[ H^1(K,E[n]) \To (R \tensor_K R)^\mult/\partial R^\mult \]
by sending $\xi \in H^1(K,E[n])$  to $\rho \, \partial R^\mult$ 
where $\rho = \partial \gamma$ for some $\gamma \in \Rbar^\mult$ such that 
$w(\xi_\sigma) = \sigma(\gamma)/\gamma$ for all $\sigma \in G_K$.

\begin{Definition} \label{DefH}
  We let $H$ denote the subgroup
  of $(R \tensor_K R)^\mult$ that maps to the image of $H^1(K, E[n])$
  in $(R \tensor_K R)^\mult/ \partial R^\mult$.
\end{Definition}

Starting with a representative $\rho \in H$, we fix a 
choice of~$\gamma$ as above. Then $\gamma$ determines a cocycle
class $\xi \in H^1(K, E[n])$. We let $\pi : C \to E$ be the twist
of the trivial $n$-covering $[n]: E \to E$ by $\xi$. 
In other words, there is a genus~1 curve $C$ defined over~$K$
and an isomorphism $\phi : C \to E$ defined over $\Kbar$
with $\sigma(\phi) \circ \phi^{-1} = \tau_{\xi_\sigma}$ (translation by
$\xi_\sigma \in E[n](\Kbar)$) for all $\sigma \in G_K$.  
The covering map is then $\pi = [n] \circ \phi$. 
It is easy to check that $\pi$ is defined over $K$.

\begin{Proposition} \label{scaling2}
  Given $\rho$ and~$C$ as above,
  there are rational functions $G_{T, C} \in \Kbar(C)$, indexed by
  $T \in E[n](\Kbar)$, such that 
  \begin{enumerate}\addtolength{\itemsep}{1mm}
    \item The divisor of $G_{T, C}$ is
           \[ \divv (G_{T, C}) = \sum_{\pi(S) = T} (S) - \sum_{\pi(S)=\oh} (S) 
                               = \pi^*(T)- \pi^*(\oh).
           \]
    \item The map $T \mapsto G_{T,C}$ is $G_K$-equivariant.
    \item The functions $G_{T, C}$ are scaled so that
          \[ r_{(T_1, T_2)}(\pi(P))
              = \rho(T_1, T_2)
                  \frac{G_{T_1, C}(P)G_{T_2, C}(P)}{G_{T_1+T_2, C}(P)} .
          \]
    \item We can package these $G_{T,C}$ to form morphisms of schemes
          \[ g_C: C \to \PP(R) \qquad \text{and} \qquad 
             g^0_C: C \setminus \pi^*E[n] \to \RR^\mult .
          \]
    \item The following diagram, with vertical maps defined over $\Kbar$, 
          commutes:
          \[ \xymatrix{ C \ar[r]^-{g_C} \ar[d]_\phi
                          & \PP(R) \ar[d]^{\cdot \gamma} \\ 
                        E \ar[r]^-{g_E} & \PP(R) } 
          \]
    \end{enumerate}
\end{Proposition}

\begin{proof}
  We take $\gamma \in \Rbar^\mult$ as above and define
  \begin{equation} \label{definproof} 
    G_{T,C}(P) = \gamma(T)^{-1} G_T(\phi(P)) 
  \end{equation}
  for $P \in C(\Kbar)$. Since $\pi =[n] \circ \phi$, statement~(i)
  is immediate from the corresponding statement for the $G_T$. 
  Next we check Galois equivariance. Since 
  $\sigma(\phi) \circ \phi^{-1} = \tau_{\xi_\sigma}$ we have 
  \[ \sigma\bigl(G_T(\phi P)\bigr) 
       = G_{\sigma T}\bigl(\phi(\sigma P) + \xi_\sigma\bigr)
       =  e_n(\xi_\sigma,\sigma T) G_{\sigma T}\bigl(\phi(\sigma P)\bigr) 
  \]
  where the second equality is the definition of
  the Weil pairing in \cite[Section~III.8]{Silverman}. 
  On the other hand, since $\sigma(\gamma)/\gamma = w(\xi_\sigma)$, we have
  \[ \sigma(\gamma(T)) = w(\xi_\sigma)(\sigma T) \gamma(\sigma T)
                       = e_n(\xi_\sigma,\sigma T) \gamma(\sigma T).
  \]
  We deduce that
  \[ \sigma(G_{T,C})(\sigma P)
       = \sigma(G_{T,C}(P)) 
       = \gamma(\sigma T)^{-1} G_{\sigma T}(\phi(\sigma P)) 
       = G_{\sigma T,C}(\sigma P) . 
  \]
  This proves~(ii). For~(iii) we compute
  \begin{align*}
    r_{(T_1,T_2)} (\pi(P)) 
      &= r_{(T_1,T_2)}(n \cdot \phi(P)) \\
      &= \frac{G_{T_1} (\phi(P)) G_{T_2} (\phi(P))}{G_{T_1+T_2}(\phi(P))} 
       = \rho(T_1, T_2) \,
           \frac{G_{T_1, C}(P) G_{T_2, C}(P)}{G_{T_1+T_2, C}(P)} .
  \end{align*}
  Here the second equality comes from Proposition~\ref{scaling}. 
  The third equality follows from~\eqref{definproof}
  and $\partial \gamma = \rho$.

  Statement~(iv) is a formal consequence of~(ii), and 
  (v) then follows from~\eqref{definproof}.
\end{proof}

The proofs of the following propositions are very similar to those of
Propositions \ref{compose} and~\ref{eqs}, and so will be omitted.

\begin{Proposition} \label{composerho}  
  The following diagram commutes.
  \[ \xymatrix{ C \setminus \pi^* E[n] \ar[rr]^{\pi} \ar[d]_{g_C^0} &
                 & E \setminus E[n] \ar[d]^{r} \\
                \RR^\mult \ar[r]^{\partial}
                  & \MS^\mult \ar[r]^{\cdot \rho} & \MS^\mult
              }
  \]
\end{Proposition}

\begin{Proposition}
  Given a Weierstrass equation for~$E$ and an element $\rho \in H$,
  with corresponding $n$-covering $\pi : C \to E$, we can explicitly compute 
  a set of $n^2 (n^2-3)/2$ linearly independent quadrics over~$K$ which
  define the image of
  \[ g_C : C \ra \PP(R) \cong \PP^{n^2-1} . \]
  If $E[n](\Kbar) = E[n](K)$, then the $z_T$ are coordinate 
  functions on $\RR$,
  and the defining quadrics can be split into two groups as follows.
  For all $T_1, T_2 \in E[n](\Kbar) \setminus \{\oh\}$, we have 
  \[ \bigl(x(T_1)-x(T_2)\bigr) z_{\oh}^2
       + \rho(T_1, -T_1) z_{T_1}z_{-T_1} - \rho(T_2, -T_2) z_{T_2}z_{-T_2} ,
  \]
  and for all $T_{11}, T_{12}, T_{21}, T_{22} \in E[n](\Kbar) \setminus \{\oh\}$
  such that 
  \[ T_{11} + T_{12} = T_{21} + T_{22} = T \neq \oh , \]
  we have
  \[ \bigl(\lambda(T_{21}, T_{22})-\lambda(T_{11}, T_{12})\bigr) z_{\oh} z_T  
       - \rho(T_{11}, T_{12}) z_{T_{11}} z_{T_{12}} 
       + \rho(T_{21}, T_{22}) z_{T_{21}} z_{T_{22}} .
  \]
\end{Proposition}


\section{A Multiplication Table for the Obstruction Algebra}\label{RcongA}

As noted in~\cite[Section~III.8]{Silverman}, for each $T \in E[n](\Kbar)$ there is a 
rational function $F_T \in \Kbar(E)$ with
\[ \divv (F_T) = n(T) - n(\oh). \]
In the last section we defined rational functions 
$G_T$. We now scale the $F_T$ so that $F_T \circ [n] = G_T^n$. 
It is equivalent to demand that the leading coefficient 
of each $F_T$, when expanded as a Laurent series in the
local parameter $x/y$ at $\oh$, should be~$1$.
(See `Step~1' in Section~5.3 of~\cite{paperI}.)  
It turns out that the $F_T$ are rather
easy to compute;  this will be explained in~\cite{paperIII}.

Following `Step~2' (loc.~cit.), we now
define $\eps \in (\Rbar \tensor_{\Kbar} \Rbar)^\mult$ by 
\[ \eps(T_1,T_2) = \frac{F_{T_1+T_2} (P)}{F_{T_1}(P)F_{T_2}(P-T_1)}. \]
for any $P \in E(\Kbar) \setminus \{\oh, T_1, T_1+T_2 \}$.
By the discussion in Section~3 of~\cite{paperI}, this does not depend on~$P$ 
and satisfies $\eps(T_1, T_2) \eps(T_2, T_1)^{-1} = e_n(T_1, T_2)$.
Since the map defining~$\eps$ is Galois-equivariant, we obtain an element
$\eps \in (R \tensor_K R)^\mult$.  The subgroup 
$H \subset (R \tensor_K R)^\mult$ was defined in Definition~\ref{DefH}.

\begin{Proposition} \label{Arho}
  There is a map
  \[ H \ra  \{\text{central simple $K$-algebras of dimension $n^2$}\} \]
  sending $\rho$ to an algebra $A_{\rho}$ such that
  \[ \Ob_n(\xi) = [A_\rho] \in \Br(K)[n], \]
  where $\xi \in H^1(K, E[n])$ is the element represented by~$\rho \in H$.
  In particular, if $\rho' = \rho \, \partial z$
  for some $z \in R^\mult$, then $[A_{\rho}] = [A_{\rho'}]$; in fact, the
  algebras $A_{\rho}$ and $A_{\rho'}$ are isomorphic.
\end{Proposition}

\begin{proof}
  In `Step~3' (loc.~cit.) we defined $A_\rho = (R,+, *_{\eps \rho})$, where
  $*_{\eps \rho}$ is a new multiplication on $R$. To define it we view
  $R \tensor_K R$ as an $R$-algebra via the comultiplication
  $\Delta : R \to R \tensor_K R$. Recall that this is defined by
  \[ \Delta(\alpha)(T_1, T_2) = \alpha(T_1 + T_2) . \]
  The corresponding trace map 
  \[ \Tr : R \tensor_K R \ra R \] 
  is given by
  $(\Tr z)(T) = \sum_{T_1+T_2 = T} z(T_1, T_2)$. Then we define
  \[ x *_{\eps \rho} y = \Tr(\eps \rho \cdot x \tensor y) \]
  for all $x, y \in R$. The stated properties of $A_\rho$ were
  established in Section~4 of~\cite{paperI}. If $\rho' = \rho\,\partial z$,
  then the isomorphism $A_{\rho'} \to A_{\rho}$ is given by
  $\alpha \mapsto z \alpha$ where the multiplication takes place in $R$.
\end{proof}

\begin{Definition}
  Let $\varphi_{\rho}: R \to A_{\rho}$ be the isomorphism of underlying
  $K$-vector spaces, inherent in the proof of Proposition~\ref{Arho}.
\end{Definition}

The construction in the proof provides us with explicit structure
constants of~$A_\rho$ in terms of a $K$-basis of~$R$.
In~\cite{paperIII}, we will discuss how to compute the structure constants
in practice.


\section{Trivialisation of the obstruction algebra} \label{BlackBox}

Recall that when $\rho\in R$ has trivial obstruction there is a
trivialisation isomorphism $\tau_\rho : A_\rho \to \Mat_n(K)$, where
$A_{\rho}$ is the central simple algebra in Proposition~\ref{Arho}.
Our algorithm will need to make this trivialisation explicit.  In
general this will be carried out by the `Black Box' to be discussed
further in~\cite{paperIII}; however, when $\rho=1$ we can write down a
standard trivialisation $\tau_1: A_1 \to \Mat_n(K)$.  In fact $\tau_1$
will depend on a choice of morphism $f_E : E \to \PP^{n-1}$ determined
by the complete linear system $|n (\oh)|$.  We make this choice now.

Let $f_E^{\vee} : E \to (\PP^{n-1})^\vee$ be the dual map of $f_E$, 
i.e., the map that takes $P \in E(\Kbar)$ to the osculating hyperplane 
at~$f_E(P)$. The elements of~$\PP^{n-1}$ will be written as column vectors,
and the elements of~$(\PP^{n-1})^\vee$ as row vectors.
For each $T \in E[n](\Kbar)$ there is a matrix $M_T \in \GL_n(\Kbar)$ such
that translation by $T$ on~$E$ extends to the automorphism of $\PP^{n-1}$ 
defined by~$M_T$. 
\begin{Proposition} \label{scaleMT}
  We may fix the scalings of the $M_T$'s so that
  \begin{equation} \label{FTandMT} 
    F_T(P) = \frac{f_E^{\vee}(\oh) \cdot M_T^{-1} \cdot f_E(P)}%
                  {f_E^{\vee}(\oh) \cdot f_E(P)} 
  \end{equation}
  where the $F_T$'s are the rational functions on~$E$ defined in 
  Section~\ref{RcongA}.  In particular, $M_{\oh} = I.$
\end{Proposition}

\begin{proof}
  The right hand side is well-defined (the undetermined scalings of
  $f_E^{\vee}(\oh)$ and of~$f_E(P)$ cancel out) and has divisor $n(T)-n(\oh)$.
\end{proof}

Let $\delta_T \in \Rbar$ be the characteristic function of $T$,
i.e., the map that takes $T$ to~$1$ but sends all other elements 
of $E[n](\Kbar)$ to~$0$. It is clear that the set of $\delta_T$ for
$T \in E[n](\Kbar)$ are a basis for $\Rbar$ as a $\Kbar$-vector space.

\begin{Definition}
  Recall the isomorphism $\varphi_1 : R \to A_1$ of underlying 
  $K$-vector spaces.
  Let $\tau_1 : A_1 \to \Mat_n(K)$ be the linear map of $K$-vector spaces 
  given by
  \[ \tau_1\bigl(\varphi_1(\alpha)\bigr)
       = \sum_{T \in E[n](\Kbar)} \alpha(T) M_T .
  \]
  (Since $T \mapsto M_T$ is $G_K$-equivariant, the map $\tau_1$,
  though {\it a priori} defined as a map of $\Kbar$-vector spaces, is
  $G_K$-equivariant.)  Note that $\tau_1$ sends $\varphi_1(\delta_T)
  \in A_1 \tensor_K \Kbar$ to $M_T \in \Mat_n(\Kbar)$.
\end{Definition}

\begin{Proposition}\label{tauisom}
  The map $\tau_1:  A_1 \to \Mat_n(K)$ is an isomorphism of $K$-algebras.
\end{Proposition}

\begin{proof}
  (See also~\cite[Prop.~5.8.(ii)]{paperI}.)
  By~\cite[Lemma 4.8]{paperI} the set of $M_T$ for $T \in E[n](\Kbar)$ 
  form a basis for $\Mat_n(\Kbar)$. So it is clear that $\tau_1$ 
  is an isomorphism of $K$-vector spaces. 
  We must show that it is also a ring homomorphism.

  We recall that $A_1 = (R,+,*_\eps)$. The new multiplication
  $*_\eps$ extends to a multiplication on~$\Rbar$ given by
  \[ \delta_{T_1} *_\eps \delta_{T_2}
       = \Tr(\eps \cdot \delta_{T_1} \tensor \delta_{T_2})
       = \eps(T_1, T_2) \delta_{T_1 + T_2} .
  \]
  Applying $\tau_1$ to both sides, it is apparent that what we have
  to show is that
  \[ M_{T_1} M_{T_2} = \eps(T_1, T_2) M_{T_1+T_2} \]
  for all $T_1, T_2 \in E[n](\Kbar)$. 
  In any case it is clear that $M_{T_1} M_{T_2} = c \, M_{T_1+T_2}$
  for some constant $c \in \Kbar^\mult$.

  Substituting $T=T_1+T_2$ in~\eqref{FTandMT}, we get
  \begin{equation} \label{FT1T2}
    F_{T_1+T_2} (P)
      = c\,\frac{f_E^{\vee}(\oh) \cdot M_{T_2}^{-1} M_{T_1}^{-1} \cdot f_E(P)}%
                {f_E^{\vee}(\oh) \cdot f_E(P)}. 
  \end{equation}
  We arbitrarily lift $f_E$ and $f^\vee_E$ to rational maps
  $E \dashrightarrow \BA^n$. The definition of $M_T$ gives
  \begin{equation} \label{defnh}
    M_{T_1}^{-1} \cdot f_E(P) = h(P) f_E(P-T_1) 
  \end{equation}
  for some rational function $h \in \Kbar(E)$. 
  Premultiplying by $f_E^\vee(\oh)$ we get
  \begin{equation} \label{evalh}
    h(P) = \frac{f_E^{\vee}(\oh) \cdot M_{T_1}^{-1} \cdot f_E(P)}%
                {f_E^\vee(\oh) \cdot f_E(P-T_1)}.
  \end{equation}
  Substituting~\eqref{defnh} into~\eqref{FT1T2} gives
  \[ F_{T_1+T_2} (P)
      = c \, \frac{f_E^{\vee}(\oh) \cdot M_{T_2}^{-1} \cdot f_E(P-T)}%
                  {f_E^{\vee}(\oh) \cdot f_E(P)} \, h(P). 
  \]
  Then by~\eqref{FTandMT} and~\eqref{evalh} we obtain
  \[ F_{T_1+T_2} (P) = c \, F_{T_2}(P-T_1) F_{T_1}(P). \]
  Comparing with the definition of~$\eps$ in Section~\ref{RcongA}, it follows 
  that $c = \eps(T_1,T_2)$ as required.
\end{proof}


\section{The Segre Factorisation for $E$}\label{SegreforE}

Recall our convention that for $V$ a vector space over~$K$, 
we write $\PP(V) = \Proj(K[V])$ where $K[V]= \oplus_{d \ge 0} \Sym^d (V^*)$
is the ring of polynomial functions on $V$.
We abbreviate $\PP(\Mat_n(K))$ as~$\PP(\Mat_n)$. The
$K$-points of~$\PP({\Mat_n})$ may be identified with the set 
$\Mat_n(K)/K^\mult$.

In the last section we fixed a morphism $f_E :E \to \PP^{n-1}$ 
determined by the complete linear system $|n (\oh)|$, and 
wrote $f_E^\vee : E \to (\PP^{n-1})^\vee$ for the dual map.
We now consider the composite map 
\[ \lambda_E : E
    \stackrel{f_E \times f_E^\vee}{\ra} \PP^{n-1} \times (\PP^{n-1})^\vee 
    \stackrel{\Segre}{\ra} \PP(\Mat_n) . 
\]

We recall that we represent elements of $\PP^{n-1}$ as column vectors,
and elements of $(\PP^{n-1})^\vee$ as row vectors.
The $\Segre$ map is then given by matrix multiplication. In particular,
its image is the locus of rank~$1$ matrices.
Since each point of $E$ lies on its own osculating hyperplane,
for $P \in E(\Kbar)$ we have $f_E^\vee(P) \cdot f_E(P) = 0$.  
Then 
\[ \Tr(\lambda_E(P)) = \Tr(f_E(P) \cdot f_E^\vee(P)) 
                     = \Tr(f_E^\vee(P) \cdot f_E(P)) = 0 .
\]
That is, the image of $\lambda_E$ is contained in the locus of trace~zero
matrices, which is a hyperplane in~$\PP(\Mat_n)$.  
There is a direct sum decomposition
$\Mat_n(K) = \langle I_n \rangle \oplus \{\Tr = 0\}$. Note that the
trace~zero subspace contains all the matrices~$M_T$ for $T \neq \oh$.
We write
$\proj$ for the rational map $\PP(\Mat_n) \dashrightarrow \PP(\Mat_n)$ 
induced by the second projection.

We defined maps $g_E : E \to \PP(R)$, $\varphi_1 : R \to A_1$ and 
$\tau_1 : A_1 \to \Mat_n(K)$ in Sections 
\ref{CtoPR}, \ref{RcongA} and \ref{BlackBox}, respectively.

\begin{Theorem} \label{segrediag} 
  The following diagram commutes.
  \[ \xymatrix{ E \ar[rr]^-{\lambda_E} \ar[d]^-{g_E} & & \PP(\Mat_n)  \\
                \PP(R) \ar[r]^-{\sim}_-{\varphi_1} 
                   & \PP(A_1) \ar[r]^-{\sim}_-{\tau_1}
                   & \PP(\Mat_n) \ar@{-->}[u]^{\operatorname{proj}} }
  \]
\end{Theorem}

The maps $g_E$ and $\tau_1$ were defined using the~$G_T$'s
and the~$M_T$'s, respectively. Note that the matrices~$M_T$ depend
on our choice of the embedding $f_E : E \to \PP^{n-1}$.

The proof of Theorem~\ref{segrediag} is based on the following result.

\begin{Theorem} \label{tomlemma} 
  For $P \in E(\Kbar) \setminus E[n](\Kbar)$ we have
  \[ \lambda_E(P) = \sum_{T \not= \oh} G_T(P) M_T. \]
\end{Theorem}

Let us show how Theorem~\ref{segrediag} follows from this.
Since the commutativity of the diagram is a geometric question, 
we are free to work over~$\Kbar$.
From the definitions we have
$\tau_1 \circ \varphi_1 \circ g_E (P) = \sum_{T} G_T(P) M_T$.
Then composing with the projection map to the trace zero subspace 
we get $P \mapsto \sum_{T \not = \oh} G_T(P) M_T$, which equals $\lambda_E(P)$
by Theorem~\ref{tomlemma}.

The proof of Theorem~\ref{tomlemma} is split into a series of lemmas.
The following notation
will be used throughout. Let $Q_1, \ldots, Q_n$ be $n$ points in~$E(\Kbar)$
and let $P =  \sum_{i=1}^n Q_i$. We define morphisms 
$h_E : E^n \to \PP^{n-1}$ and $h_E^\vee : E^n \to (\PP^{n-1})^\vee$
as follows. For the first map we put
\[ h_E(Q_1, \ldots, Q_n) = f_E(P). \]
The second map takes $(Q_1, \dots, Q_n)$ to the
hyperplane meeting $E$ (or rather the image of $f_E$) 
in the divisor $(P - n Q_1) + \ldots + (P - n Q_n)$. Notice that
since this divisor has sum~$\oh$, it is indeed linearly
equivalent to the hyperplane section.

The first lemma can be seen as a multi-variable generalisation of 
formula~\eqref{FTandMT} in Proposition~\ref{scaleMT}.

\begin{Lemma} \label{step1lemma}
  If $Q_1, \ldots, Q_n \in E(\Kbar) \setminus E[n](\Kbar)$ then  
  \[ G_T(Q_1) \ldots G_T(Q_n)
      = \frac{h_E^\vee(Q_1,\ldots,Q_n)\cdot M_T^{-1}\cdot h_E(Q_1,\ldots,Q_n)}%
             {h_E^\vee(Q_1,\ldots,Q_n) \cdot h_E(Q_1,\ldots,Q_n)} 
  \]
  for all~$T \in E[n](\Kbar)$. 
\end{Lemma}

\begin{proof}
  We view each side as a rational function on $E^n$.
  The strategy of the proof is first to compare divisors, and
  then to check scalings by specialising to the case
  $Q_1 = Q_2 = \ldots = Q_n$. 

  Let $\pr_i : E^n \to E$ be projection to the $i$th factor.
  The left hand side has divisor 
  \[ \sum_{i=1}^n \sum_{nx =T} \pr_i^*(x) 
       - \sum_{i=1}^n \sum_{nx=\oh} \pr_i^* (x). 
  \]
  From the definitions of $h_E$ and~$h_E^\vee$ the right hand side 
  has a zero whenever $n Q_i = T$ for some~$i$, and a pole whenever 
  $n Q_i = \oh$ for some~$i$. Therefore the right hand 
  side has divisor
  \[ \sum_{i=1}^n \sum_{nx=T} a_{x} \pr_i^*(x) 
       - \sum_{i=1}^n  \sum_{nx=\oh} b_{x} \pr_i^* (x) 
  \]
  where the $a_{x}$ and $b_{x}$ are positive integers.

  If we replace $Q_1$ by $Q_1 + S$ for $S \in E[n](\Kbar)$,
  then the right hand side is multiplied by a non-zero scalar
  (the commutator of $M_{S}$ and $M_T$).
  It follows that the integers $a_{x}$ and $b_{x}$ do not depend on~$x$. 
  So the right hand side has divisor
  \[ a \sum_{i=1}^n \sum_{nx=T} \pr_i^*(x) 
       - b \sum_{i=1}^n \sum_{nx=\oh} \pr_i^* (x) 
  \]
  where $a$ and $b$ are positive integers. 

  Since this divisor is principal, its pull-back by any morphism
  $E \to E^n$ has degree $0$.
  This enables us to show that $a = b$. 
  Since $E^n$ is a projective variety, 
  it follows that the right hand side is
  \[ c \, G_T(Q_1)^{a} \cdots G_T(Q_n)^{a} \]
  for some constant $c \in \Kbar^\mult$. 
  
  We now specialise by taking $Q_1= Q_2 = \ldots = Q_n$ ($= Q$ say).
  Then $P=nQ$ and 
  \[ c \, G_T(Q)^{n a}
       =  \frac{f_E^\vee(\oh) \cdot M_T^{-1} \cdot f_E(P)}%
               {f_E^\vee(\oh) \cdot f_E(P)} . 
  \]
  The definition of~$F_T$ in Section~\ref{RcongA} gives 
  $F_T(P)= F_T(nQ) = G_T(Q)^n$.
  Finally we compare with~\eqref{FTandMT} to get $c = 1$ and $a = 1$.
\end{proof}

\begin{Lemma} \label{step2lemma}
  If $Q_1, \ldots, Q_n \in E(\Kbar) \setminus E[n](\Kbar)$ then 
  \[ \frac{1}{n} \sum_{T \in E[n](\Kbar)} G_T (Q_1) \ldots G_T(Q_n) M_T 
       = \frac{h_E(Q_1, \ldots, Q_n) \cdot h^\vee_E(Q_1, \ldots, Q_n)}%
              {h_E^\vee(Q_1, \ldots, Q_n) \cdot h_E(Q_1, \ldots, Q_n)}
  \]
\end{Lemma}

\begin{proof}
  Since the $M_T$ form a basis for $\Mat_n(\Kbar)$, we can write
  the right hand side as $n^{-1} \sum_T a_T(Q_1, \ldots, Q_n) M_T$ 
  for some rational functions $a_T$ on~$E^n$. To compute the
  $a_T$'s we premultiply by $M_T^{-1}$ and take the trace.
  Since 
  \[ \Tr(M_T) = \begin{cases}
                  n & \text{ if $T = \oh$,} \\
                  0 & \text{ otherwise,} 
                \end{cases}
  \]
  on the left hand side we get $G_T (Q_1) \ldots G_T(Q_n)$.
  On the right hand side, using $\Tr(AB) = \Tr(BA)$ gives
  \begin{align*}
    a_T(Q_1, \ldots, Q_n)
      &= \Tr\Bigl(M_T^{-1} \cdot
                  \frac{h_E(Q_1, \ldots, Q_n) \cdot h^\vee_E(Q_1, \ldots, Q_n)}%
                       {h_E^\vee(Q_1, \ldots, Q_n) \cdot h_E(Q_1, \ldots, Q_n)}
            \Bigr) \\
      &= \Tr\Bigl(\frac{h_E^{\vee}(Q_1, \ldots, Q_n) \cdot M_T^{-1}
                                              \cdot h_E(Q_1, \ldots, Q_n)}%
                       {h_E^\vee(Q_1, \ldots, Q_n) \cdot h_E(Q_1, \ldots, Q_n)}
            \Bigr) \\
      &= G_T (Q_1) \ldots G_T(Q_n) ,
  \end{align*}
  where we have used Lemma~\ref{step1lemma} in the last equality.
\end{proof}

\begin{Lemma} \label{step3lemma}
  There is a commutative diagram of morphisms
  \[ \xymatrix{ E^n \ar[r]^-{h_E \times h_E^\vee} \ar[d]^{\prod g_E}
                  & \PP^{n-1} \times (\PP^{n-1})^\vee \ar[r]^-{\Segre} 
                  & \PP(\Mat_n) \\ 
                \PP(R) \ar[r]^-{\sim}_-{\varphi_1}
                  & \PP(A_1) \ar[r]^-{\sim}_-{\tau_1} 
                  & \PP(\Mat_n) \ar@{=}[u]  }
  \]
  where $(\prod g_E) (Q_1, \ldots, Q_n)= \prod_{i=1}^n g_E(Q_i)$.
\end{Lemma}

\begin{proof}
  From the definitions we have
  \[ \tau_1 \circ \varphi_1 \circ (\prod g_E) (Q_1, \ldots ,Q_n)
       = \sum_{T} G_T(Q_1) \ldots G_T(Q_n) M_T .
  \]
  So the commutativity is already clear from Lemma~\ref{step2lemma}.

  It only remains to check that $\prod g_E$ is a morphism. To do
  this we write it as a composite 
  \[ E^n \stackrel{g_E^n}{\ra} \PP(R)^n  
         \stackrel{\scriptscriptstyle\prod}{\dashrightarrow} \PP(R) \]
  where the second map is induced by multiplication in~$R$. We check that 
  the image of the first map is contained in the domain 
  of definition of the second. 

  If $Q \in E(\Kbar) \setminus E[n](\Kbar)$, then 
  $g_E(Q)$ is the class of $T \mapsto G_T(Q)$,
  whereas if $Q \in E[n](\Kbar)$, then $g_E(Q)$ is the 
  class of $T \mapsto \res_Q(G_T)$, where the residue
  is taken with respect to a local parameter at~$Q$. (The choice of local
  parameter does not matter.) 

  Now suppose $Q_1, \ldots, Q_n \in E(\Kbar)$ and $\prod_{i=1}^n g_E(Q_i)$
  is undefined as an element of~$\PP(R)$. Then for each $T \neq \oh$ there
  exists $1 \le i \le n$ such that $G_T(Q_i) = 0$, and hence $n Q_i = T$. 
  But there
  are $n^2 - 1$ such choices of $T$ and only $n$ choices of~$i$. 
  So this is impossible. It follows that $\prod g_E$ is a morphism as claimed.
\end{proof}

To complete the proof of Theorem~\ref{tomlemma}, 
and hence of Theorem~\ref{segrediag}, 
we put $Q_1= Q$ and $Q_2 = \ldots = Q_n = \oh$ in Lemma~\ref{step3lemma}.
Notice that $h_E(Q, \oh, \ldots, \oh) = f_E(Q)$
and $h_E^\vee(Q, \oh, \ldots, \oh) = f_E^\vee(Q)$. 
Also, with the $G_T$'s scaled as in Proposition~\ref{scaling}(iii),
$g_E(\oh)$ is the class of $\sum _{T \neq \oh} \delta_T$.
So for $Q \in E(\Kbar) \setminus E[n](\Kbar)$ we obtain
\begin{align*}
  \lambda_E(Q)
    &= {\Segre} \circ (f_E \times f_E^{\vee})(Q) \\
    &= {\Segre} \circ (h_E \times h_E^{\vee})(Q, \oh, \dots, \oh) \\
    &= \tau_1 \circ \varphi_1 (g_E(\oh)^{n-1} g_E(Q)) \\
    &= \sum_{T \not= \oh} G_T(Q) M_T 
\end{align*}


\section{The Segre Factorisation for $C$} \label{SegreforC}

Let $A$ be a central simple algebra over~$K$ of dimension~$n^2$.
Let $S$ and~$S^\vee$ be the Brauer-Severi varieties given by
the minimal right and left ideals of $A$, respectively 
(see~\cite[p.~160]{Serre}). There is a natural map
$\Segre : S \times S^\vee \to \PP(A)$ given by intersecting ideals. 
We say an element $a \in A$ has rank~$r$ if the map $x \mapsto a x$
is an endomorphism of~$A$ (as a~$K$-vector space) of rank~$rn$.

\begin{Lemma} \label{whyrank1good}
  The $\Segre$ map $S \times S^\vee \to \PP(A)$ is an embedding,
with image the locus of rank~1 elements in $\PP(A)$.
\end{Lemma}

\begin{proof} 
  If $A \isom \Mat_n(K)$ then $S \isom \PP^{n-1}$ and the Segre map
  reduces to that studied in Section~\ref{SegreforE}. 
  The description of the image
  is no more than the observation that a (non-zero) matrix has rank~1 
  if and only if it can be written as a column vector times a row vector. 
  In general we use that there is an isomorphism of $\Kbar$-algebras
  $A \tensor_K \Kbar \isom \Mat_n(\Kbar)$. 

  The Segre map is an embedding since, on the rank~1 locus in $\PP(A)$,
  an inverse is given by sending the class of~$a \in A$ to the pair 
  of ideals $(aA, Aa)$.
\end{proof}

From now on, we call the Segre map the (generalised) Segre embedding.
We write $1_A$ for the multiplicative identity of~$A$, and 
$\Trd : A \to K$ for the reduced trace.
There is a decomposition of $K$-vector spaces 
$A = \langle 1_A \rangle \oplus \{ \Trd = 0 \}$.
As in Section~\ref{SegreforE} we write $\proj$ 
for projection onto the second factor.

The subgroup $H \subset (R \tensor_K R)^\mult$ was defined in 
Definition~\ref{DefH}.

\begin{Theorem} \label{twistedsegrediag} 
  Let $\rho \in H$, and let $\pi :C \to E$ be the corresponding $n$-covering.
  Let $S$ and $S^\vee$ be the Brauer-Severi varieties given by 
  the minimal right and left ideals in~$\PP(A_\rho)$.
  Then there is a morphism $f_C : C \to S$ with dual
  $f_C^\vee : C \to S^\vee$ such that 
  \begin{enumerate}
    \item the following diagram commutes:
          \[ \xymatrix{ C \ar[r]^-{f_C \times f_C^\vee} \ar[d]^{g_C} 
                          & S \times S^\vee \ar[r]^-{\Segre} 
                          & \PP(A_{\rho}) \\ 
                        \PP(R) \ar[rr]^-{\sim}_-{\varphi_\rho}
                          & & \PP(A_{\rho}) \ar@{-->}[u]^{\operatorname{proj}}
                      } 
          \]
    \item the Brauer-Severi diagram $[C \to S]$ and the class of $\rho$ in $H$
          correspond to the same element of $H^1(K, E[n])$.
  \end{enumerate}
\end{Theorem}

The theorem is proved by combining Theorem~\ref{segrediag} with the 
next lemma. First we recall how the element $\rho \in H$ and
$n$-covering $\pi :C \to E$ are related. In Section~\ref{CtoPR} 
we fixed an element $\gamma \in \Rbar^\mult$ with $\rho = \partial \gamma$
and defined a cocycle $\xi \in H^1(K,E[n])$ via 
$w(\xi_\sigma) = \sigma(\gamma)/\gamma$ for all $\sigma \in G_K$.
Then we let $\pi : C \to E$ be the $\xi$-twist
of the trivial $n$-covering. 
Thus there is an isomorphism $\phi : C \to E$ defined over $\Kbar$
with $\pi = [n] \circ \phi$ and 
$\sigma(\phi) \circ \phi^{-1} = \tau_{\xi_\sigma}$ 
for all $\sigma \in G_K$.

\begin{Lemma} \label{landinrank1}
There is an isomorphism of $\Kbar$-algebras 
$\beta : A_\rho \otimes_K \Kbar \to A_1 \otimes_K \Kbar$ making the
following diagram commute.
\begin{equation} \label{newdiagram}
       \xymatrix{ C \ar[r]^-{g_C} \ar[d]^{\phi}
                  & \PP(R) \ar[r]^{\sim}_-{\varphi_\rho} \ar[d]^{\cdot \gamma} 
                  & \PP(A_\rho) \ar[d]^{\beta} 
                    \ar[dr]^{\tau_1 \circ \beta} \ar@{-->}[r]^-{\proj} 
                  & \PP(A_\rho) \ar[dr]^{\tau_1 \circ \beta} \\
                E \ar[r]^-{g_E}
                  & \PP(R) \ar[r]^{\sim}_-{\varphi_1} 
                  & \PP(A_1) \ar[r]^-{\tau_1} 
                  & \PP(\Mat_n) \ar@{-->}[r]^-{\proj} 
                  & \PP(\Mat_n) }
\end{equation}
\end{Lemma}

\begin{proof} 
  The first square commutes by Proposition~\ref{scaling2}(v). 
  We define $\beta : A_{\rho} \tensor_K \Kbar \to A_1 \tensor_K \Kbar$ 
  to make the second square commute. 
  It is an isomorphism of $\Kbar$-algebras by~\cite[Lemma 4.6]{paperI}. 
  We recall from Proposition~\ref{tauisom} that $\tau_1$ is
  an isomorphism of $K$-algebras. Finally, since $\proj$ is defined 
  purely in terms of the algebra structure, it commutes with 
  the algebra isomorphism $\tau_1 \circ \beta$.
\end{proof}

Theorem~\ref{segrediag} identifies the composite of the second row 
of~(\ref{newdiagram}) as $\lambda_E$, where we recall
\[ \lambda_E : E
    \stackrel{f_E \times f_E^\vee}{\ra} \PP^{n-1} \times (\PP^{n-1})^\vee 
    \stackrel{\Segre}{\ra} \PP(\Mat_n) . \]
Since $\lambda_E$ factors via the Segre embedding, its image
belongs to the rank~1 locus of $\PP(\Mat_n)$. Let 
$\lambda_C = \proj \circ \varphi_{\rho} \circ g_C$ be the
composite of the first row of~(\ref{newdiagram}). Recalling that
$\tau_1 \circ \beta$ is an isomorphism of $\Kbar$-algebras,
it follows that $\lambda_C$ has image belonging to the rank~1
locus of $\PP(A_\rho)$. By Lemma~\ref{whyrank1good} it therefore
factors via the Segre embedding.

In other words, there are morphisms
$f_C : C \to S$ and $f_C^{\vee} : C \to S^{\vee}$ making the diagram
in the first part of Theorem~\ref{twistedsegrediag} commute.
It remains to show that $[C \to S]$ is a Brauer-Severi diagram,
that $f_C^{\vee}$ is dual to~$f_C$ and that $[C \to S]$ and~$\rho$
correspond to the same element of~$H^1(K, E[n])$.

The isomorphism of $\Kbar$-algebras
$\tau_1 \circ \beta : A_\rho \tensor_K \Kbar \to \Mat_n(\Kbar)$
induces isomorphisms $\psi : S \to \PP^{n-1}$ and
$\psi^{\vee} : S^{\vee} \to (\PP^{n-1})^{\vee}$ defined over~$\Kbar$
making the following diagram commute.
\[ \xymatrix{ C \ar[rr]^-{f_C \times f_C^{\vee}} \ar[d]_{\phi}
                & & S \times S^{\vee} \ar[rr]^-{\Segre}
                                      \ar[d]_{\psi \times \psi^{\vee}}
                & & \PP(A_\rho) \ar[d]^{\tau_1 \circ \beta} \\
              E \ar[rr]^-{f_E \times f_E^{\vee}}
                & & \PP^{n-1} \times (\PP^{n-1})^{\vee} \ar[rr]^-{\Segre}
                & & \PP(\Mat_n)
            }
\]
This shows that via $(\phi, \psi)$, the morphism $f_C : C \to S$ is isomorphic
over~$\Kbar$ to $f_E : E \to \PP^{n-1}$.  Hence $[C \to S]$
is a Brauer-Severi diagram.  Since $\sigma(\phi) \circ \phi^{-1}$
is translation by $\xi_\sigma$ this diagram is 
the $\xi$-twist of $[E \to \PP^{n-1}]$,
and therefore the Brauer-Severi diagram corresponding to~$\rho$.
At the same time, we see that $f_C^{\vee}$
is dual to~$f_C$ (since $f_E^{\vee}$ is dual to~$f_E$).
This completes the proof of Theorem~\ref{twistedsegrediag}.


\section{Finding Equations for $C$ in $\PP^{n-1}$} \label{CinPn}

We continue to represent elements of $H^1(K,E[n])$
by elements $\rho \in H$, where the subgroup 
$H \subset (R \otimes_K R)^\mult$ was defined in Section~\ref{CtoPR}.
The obstruction algebra $A_\rho$ was introduced 
in Section~\ref{RcongA}.
If $\rho$ represents an element of $H^1(K,E[n])$ with trivial
obstruction, then we make use of an explicit trivialisation
isomorphism of $K$-algebras $\tau_{\rho} : A_\rho \to \Mat_n(K)$.

\begin{Theorem} \label{masterdiagram}
  Let $\rho \in H$, and let $\pi :C \to E$ be the corresponding $n$-covering.
  Suppose given an isomorphism of $K$-algebras 
  $\tau_{\rho} : A_\rho \to \Mat_n(K)$.
  Then there is a morphism $f'_C : C \to \PP^{n-1}$ with dual ${f'_C}^{\vee}$
  such that 
  \begin{enumerate}
    \item the following diagram commutes:
   \[ \xymatrix{ C \ar[r]^-{f'_C \times {f'_C}^\vee} \ar[d]^{g_C}
                  & \PP^{n-1} \times (\PP^{n-1})^\vee \ar[r]^-{\Segre} 
                  & \PP(\Mat_n) \\
                \PP(R) \ar[r]_-{\varphi_\rho}^-{\sim}
                  & \PP(A_\rho) \ar[r]_-{\tau_\rho}^-{\sim}
                  & \PP(\Mat_n) \ar@{-->}[u]^{\proj} }
  \]
    \item the Brauer-Severi diagram $[C \to \PP^{n-1}]$ and the class 
    of $\rho$ in $H$ correspond to the same element of $H^1(K, E[n])$.
  \end{enumerate}
\end{Theorem}

\begin{proof} 
  Let $S$ and $S^\vee$ be the Brauer-Severi varieties given by 
  the minimal right and left ideals in~$\PP(A_\rho)$.
  The isomorphism $\tau_\rho$ induces $K$-isomorphisms
  $\psi' : S \to \PP^{n-1}$ and 
  ${\psi'}^{\vee} : S^{\vee} \to (\PP^{n-1})^{\vee}$.
  We modify the maps $f_C$ and $f_C^{\vee}$ of 
  Theorem~\ref{twistedsegrediag} to give maps $f'_C$ and ${f'_C}^{\vee}$ 
  making the following diagram commute.
  \[ \xymatrix{ C \ar[rr]^-{f'_C \times {f'_C}^{\vee}}
                  & & \PP^{n-1} \times (\PP^{n-1})^{\vee} \ar[rr]^-{\Segre}
                  & & \PP(\Mat_n) \\
                C \ar[rr]^-{f_C \times f_C^{\vee}} \ar@{=}[u]
                  & & S \times S^{\vee} \ar[rr]^-{\Segre}
                                        \ar[u]_{\psi' \times {\psi'}^{\vee}}
                  & & \PP(A_\rho) \ar[u]_{\tau_\rho}
              }
  \]
  Since $\proj$ is defined purely in terms of the algebra
  structure, it commutes with the algebra isomorphism $\tau_{\rho}$.
  The theorem now follows by combining the above diagram with that in
  Theorem~\ref{twistedsegrediag}.
\end{proof}

In Section~\ref{CtoPR} we explained how to write down equations
for $C$ as a curve of degree $n^2$ in $\PP(R) \isom \PP^{n^2-1}$.
Theorem~\ref{masterdiagram} now tells us how to convert these
to equations for $C$ as a curve of degree $n$ in $\PP^{n-1}$. 
First we use $\tau_{\rho}$ to get equations for $C$ in~$\PP(\Mat_n)$.
Then we project onto the trace zero matrices. 
Next we write $x_{11}, x_{12}, \ldots, x_{nn}$ for our coordinate
functions on $\PP(\Mat_n)$ and substitute $x_{ij} = x_i y_j$, 
where the $x_i$ and $y_j$ are new indeterminates ($2n$ in total). 
This corresponds to pulling the image of~$C$ in~$\PP(\Mat_n)$ 
back under the Segre map.
To project to the first factor~$\PP^{n-1}$, we eliminate the~$y_j$.
We are left with equations in the~$x_i$, and these define the 
image of $f'_C : C \to \PP^{n-1}$. We will describe in~\cite{paperIII}
how this computation can be reduced to linear algebra 
for any specific value of~$n$.

For aesthetic as well as practical reasons, it is desirable to 
find a change of coordinates on $\PP^{n-1}$ so 
that the equations for $C$ have small
coefficients. The necessary minimisation
and reduction procedures will be described in a forthcoming 
paper~\cite{paperIV} for the cases $n = 3$ and~$n = 4$.

It is also desirable to have explicit equations for the 
covering map $\pi : C \to E$. In principle these could be obtained
using the methods of Section~\ref{CtoPR}. However in
the cases $n=2,3,4$ equations for $\pi$ are already given by 
classical formulae, reproduced in \cite{AKM3P}.


\end{document}